\documentclass{amsart}

\usepackage{amssymb,amsmath,epsfig}
\usepackage{amsfonts}

\newtheorem{theorem}{Theorem}[section]
\newtheorem{lemma}[theorem]{Lemma}

\theoremstyle{definition}

\theoremstyle{remark}

\numberwithin{equation}{section}

\begin{document}

\title{ Two identities involving the cubic partition function}

\author{Xinhua, xiong}
\address{Department of Mathematics, China Three Gorges University, Yichang, Hubei Province, 443002, P.R. China}

\email{xinhuaxiong@ctgu.edu.cn}


\date{}

\commby{}

\begin{abstract}
 We give a new proof of
Chan's identity involving the cubic partition function and we also give a new identity for the cubic partition function which is analogues to the Zuckerman's identity for the ordinary partition function.
\end{abstract}

\maketitle

\section{INTRODUCTION}
Let $p(n)$ be the number of partitions of $n$, defined by $\sum_{n\geq 0}p(n)q^n:= \prod_{n=1}^{\infty}(1-q^n)^{-1}$.
In connection with his discovery of certain divisibility properties of $p(n)$ Ramanujan  stated the identities:
\begin{equation}\label{ramanujan1}
\sum_{n=0}^{\infty}p(5n+4)q^n = 5\frac{(q^5;q^5)_{\infty}^5}{(q;q)_{\infty}^6},
\end{equation}
and
\begin{equation}\label{ramanujan2}
\sum_{n=0}^{\infty}p(7n+5)q^n = 7\frac{(q^7;q^7)_{\infty}^3}{(q;q)_{\infty}^4}+49q\frac{(q^7;q^7)_{\infty}^7}{(q;q)_{\infty}^8}.
\end{equation}
Here and in the rest of the paper we follow the customary $q$-product notation: we set (for $|q| \le1$)
$$
(a;q)_{\infty}= \prod_{n=0}^{\infty}(1-aq^n).
$$
Both Hardy and MacMahon considered \ref{ramanujan1} as Ramanujan's "Most Beautiful Identity". Darling \cite{Dar21} proved the first and Mordell \cite{Mo22}, Watson \cite{Watson38} and Rademacher-Zuckerman \cite{Ra38} gave proofs for both identities. Recently, H.H. Chan and R.P. Lewis \cite{H-Chan} also gave different proofs of both identities. All these proofs used the theory of modular functions. In another paper, Zuckerman \cite{Zu41} obtained the following identity (Zuckerman's identity):
\begin{eqnarray}\label{ramanujan3}
\sum_{n=0}^{\infty}p(25n+24)q^n = 63\cdot 5^2 \frac{(q^5;q^5)_{\infty}^6}{(q;q)_{\infty}^7}&+&52\cdot 5^5 q \frac{(q^5;q^5)_{\infty}^{12}}{(q;q)_{\infty}^{13}}+63\cdot 5^7 q^2 \frac{(q^5;q^5)_{\infty}^{18}}{(q;q)_{\infty}^{19}}
\cr&+&6\cdot 5^{10} q^3 \frac{(q^5;q^5)_{\infty}^{24}}{(q;q)_{\infty}^{25}}+5^{12}q^4 \frac{(q^5;q^5)_{\infty}^{30}}{(q;q)_{\infty}^{31}}.
\end{eqnarray}

In  two recent papers,
H.-C. Chan \cite{Chan08a, Chan10} proved a generalization of \ref{ramanujan1} and \ref{ramanujan2} for a certain kind of partition function $a(n)$ which is defined by
\begin{equation}
\sum_{n=0}^{\infty}a(n)q^n:= \frac{1}{(q;q)_{\infty}(q^2;q^2)_{\infty}}.
\end{equation}
Kim \cite{Kim08} noted that $a(n)$ can be interpreted the number of $2$-color partitions of $n$ with colors r and g subject to the restriction that the color b appears only in even parts,
so he called $a(n)$ to be the cubic partition function owing to the fact that $a(n)$ is related to Ramanujan's cubic
continued fraction. Using some identities for the cubic continued fraction, H.C.-Chan derived the following identity :
\begin{theorem}[\cite{Chan08a} Theorem 1]
\begin{equation}\label{chanxizhi}
\sum_{n=0}^{\infty}a(3n+2)q^n=\frac{3(q^3;q^3)_\infty^3(q^6;q^6)_\infty^3}{(q;q)_{\infty}^4(q^2;q^2)_{\infty}^4}.
\end{equation}
\end{theorem}
In this note, we only use 3-dissections of functions $ \frac{1}{\Phi(-q)} $ and $ \frac{1}{\Psi(q)} $ to give a new elementary proof
of \ref{chanxizhi} and we also give the following identity for the cubic partition function $a(n)$ which is similar to  Zuckerman's identity \ref{ramanujan3} for the cubic partition function $a(n)$.
\begin{theorem}
\begin{eqnarray}\label{xiong}
\sum_{n=0}^{\infty}a(9n+8)q^n &=& 2\cdot 3^3 \frac{(q^3;q^3)_{\infty}^{30}}{(q;q)_{\infty}^{19}(q^2;q^2)_{\infty}^{7}(q^6;q^6)_{\infty}^{6}}+8\cdot 3^3 q \frac{(q^3;q^3)_{\infty}^{21}(q^6;q^6)_{\infty}^{3}}{(q;q)_{\infty}^{16}(q^2;q^2)_{\infty}^{10}}\cr\cr&+&19\cdot 3^4 q^2 \frac{(q^3;q^3)_{\infty}^{12}(q^6;q^6)_{\infty}^{12}}{(q;q)_{\infty}^{13}(q^2;q^2)_{\infty}^{13}}
-64\cdot 3^{3} q^3 \frac{(q^3;q^3)_{\infty}^{3}(q^6;q^6)_{\infty}^{21}}{(q;q)_{\infty}^{10}(q^2;q^2)_{\infty}^{16}}\cr\cr&+&128\cdot 3^{3}q^4 \frac{(q^6;q^6)_{\infty}^{30}}{(q;q)_{\infty}^{7}(q^2;q^2)_{\infty}^{19}(q^3;q^3)_{\infty}^{6}}.
\end{eqnarray}
\end{theorem}
\noindent{Hence $a(9n+8)\equiv 0\pmod{27}$, which coincides with the result of Chan, he derived this result with different method.}

\section{PRELIMINARIES}
We require a few definitions and lemmas. Let
$$
\Phi(-q)=\sum_{n=-\infty}^{\infty}(-1)^nq^{n^2}=\frac{(q;q)_{\infty}^{2}}{(q^2;q^2)_{\infty}},\quad \Psi(q)=\sum_{n=0}^{\infty}q^{(n^2+n)/2}=\frac{(q^2;q^2)_{\infty}}{(q;q^2)_{\infty}},
$$
and
$$P(q)=\frac{(q^2;q^6)_\infty(q^4;q^6)_\infty(q^3;q^3)_{\infty}^2}{(q;q)_\infty},\quad
X(-q)=\frac{(q;q)_{\infty}^{}(q^6;q^6)_{\infty}^{2}}{(q^2;q^2)_{\infty}^{}(q^3;q^3)_{\infty}^{}}.
$$
We will use the following necessaries in proving our main results.
The following lemma are the 3-dissections of functions $ \frac{1}{\Phi(-q)} $ and $ \frac{1}{\Psi(q)} $.
\begin{lemma}[\cite{Hir1}, last line in the proof Theorem 1]\label{lemma2.2}
$$
\frac{1}{\Phi(-q)}=\frac{\Phi(-q^9)}{{\Phi(-q^3)}^4}(\Phi(-q^9)^2+2q \Phi(-q^9)X(-q^3)+4q^2 X(-q^3)^2).
$$
\end{lemma}
\begin{lemma}[\cite{Hir2}, Lemma 2.2]\label{lemma2.3}
$$
\frac{1}{\Psi(q)}=\frac{\Psi(q^9)}{{\Psi(q^3)}^4}(P(q^3)^2-q P(q^3)\Psi(q^9)+q^2 \Psi(q^9)^2).
$$
\end{lemma}
\begin{lemma}\label{lemma2.4}
$$
\Phi(-q)\Psi(q)=(q;q)_{\infty}(q^2;q^2)_{\infty},\quad
X(-q)P(q)=(q^3;q^3)_{\infty}(q^6;q^6)_{\infty}.
$$
\end{lemma}
\begin{proof}
We have
\begin{eqnarray}\label{2.1}
\Phi(-q)\Psi(q)=\frac{(q;q)_{\infty}^2}{(q^2;q^2)_{\infty}}\cdot\frac{(q^2;q^2)_{\infty}}
{(q;q^2)_{\infty}}=\frac{(q;q)_{\infty}^2}{(q;q^2)_{\infty}}=(q;q)_{\infty}(q^2;q^2)_{\infty}.
\end{eqnarray}
and
\begin{eqnarray}\label{2.2}
X(-q)P(q)&=&\frac{(q;q)_{\infty}^{}(q^6;q^6)_{\infty}^{2}}{(q^2;q^2)_{\infty}^{}(q^3;q^3)_{\infty}^{}}
\cdot \frac{(q^2;q^6)_\infty(q^4;q^6)_\infty(q^3;q^3)_{\infty}^2}{(q;q)_\infty}\cr\cr
&=& \frac{(q^2;q^6)_{\infty}(q^4;q^6)_{\infty}(q^3;q^3)_{\infty}(q^6;q^6)_{\infty}^2}{(q^2;q^2)_{\infty}}\cr\cr
&=& \frac{(q^2;q^6)_{\infty}(q^4;q^6)_{\infty}(q^3;q^3)_{\infty}(q^6;q^6)_{\infty}^2}
{(q^2;q^6)_{\infty}(q^4;q^6)_{\infty}(q^6;q^6)_{\infty}}\cr\cr &=& (q^3;q^3)_{\infty}(q^6;q^6)_{\infty}.
\end{eqnarray}
\end{proof}

\section{ AN  NEW PROOF OF IDENTITY \ref{chanxizhi}  }
We begin with the proof of identity \ref{chanxizhi}.
\begin{proof}
We note that
$$
\frac{1}{(q;q)_{\infty}(q^2;q^2)_\infty}=\frac{(q^2;q^2)_{\infty}}{(q;q)_{\infty}^2}\cdot \frac{(q;q^2)_{\infty}}{(q^2;q^2)_{\infty}}=\frac{1}{\Phi(-q)\Psi(q)}
$$
by using lemma \ref{lemma2.3}. So we have
\begin{eqnarray*}
&&\sum_{n=0}^\infty a(n)q^n =\frac{1}{\Phi(-q)\Psi(q)}\cr
&=& \frac{\Phi(-q^9)}{{\Phi(-q^3)}^4}(\Phi(-q^9)^2+2q \Phi(-q^9)X(-q^3)+4q^2 X(-q^3)^2)\cr &&
\frac{\Psi(q^9)}{{\Psi(q^3)}^4}(P(q^3)^2-q P(q^3)\Psi(q^9)+q^2 \Psi(q^9)^2)\cr
&=& \frac{\Phi(-q^9)\Psi(q^9)}{\Phi(-q^3)^4\Psi(q^3)^4}(\Phi(-q^9)^2P(q^3)^2+2q^3 \Phi(-q^9)\Psi(q^9)^2X(-q^3)\cr
&-&4q^3\Psi(q^9)X(-q^3)^2P(q^3)
+ 2q\Phi(-q^9)X(-q^3)P(q^3)^2-q \Phi(-q^9)^2\Psi(q^9)P(q^3)\cr
&+&4q^4 \Psi(q^9)^2 X(-q^3)^2
+ q^2 \Phi(-q^9)^2 \Psi(q^9)^2+ 4q^2X(-q^3)^2 P(q^3)^2\cr
&-& 2q^2 \Phi(-q^9)\Psi(q^9)X(-q^3)P(q^3)).
\end{eqnarray*}
Therefore,
\begin{eqnarray*}
&&\sum_{n=0}^\infty a(3n+2)q^{3n+2}\cr&=&\frac{\Phi(-q^9)\Psi(q^9)}{\Phi(-q^3)^4\Psi(q^3)^4}(q^2 \Phi(-q^9)^2 \Psi(q^9)^2+
 4q^2X(-q^3)^2 P(q^3)^2\cr
 &-& 2q^2 \Phi(-q^9)\Psi(q^9)X(-q^3)P(q^3)).
\end{eqnarray*}
Dividing by $q^2$ on both sides and replacing $q^3$ by $q$, we obtain
\begin{eqnarray*}
&&\sum_{n=0}^\infty a(3n+2)q^{n}\cr&=&\frac{\Phi(-q^3)\Psi(q^3)}{\Phi(-q)^4\Psi(q)^4}( \Phi(-q^3)^2 \Psi(q^3)^2+
 4X(-q)^2 P(q)^2- 2 \Phi(-q^3)\Psi(q^3)X(-q)P(q))\cr
 &=& \frac{(q^3;q^3)_\infty^3(q^6;q^6)_\infty^3}{(q;q)_{\infty}^4(q^2;q^2)_{\infty}^4} +4\frac{(q^3;q^3)_\infty(q^6;q^6)_\infty}{(q;q)_{\infty}^4(q^2;q^2)_{\infty}^4}(q^3;q^3)_\infty^2(q^6;q^6)_\infty^2\cr
 &-&2\frac{(q^3;q^3)_\infty^2(q^6;q^6)_\infty^2}{(q;q)_{\infty}^4(q^2;q^2)_{\infty}^4}(q^3;q^3)_\infty(q^6;q^6)_\infty\cr
 &=& 3\frac{(q^3;q^3)_\infty^3(q^6;q^6)_\infty^3}{(q;q)_{\infty}^4(q^2;q^2)_{\infty}^4},
 \end{eqnarray*}
by using Lemma \ref{lemma2.4}.
\end{proof}
\vskip -1cm
\section{PROOF OF THE IDENTITY \ref{xiong}}
We will use the following lemma in proving the identity \ref{xiong}.
\begin{lemma}\label{4.1}
Let $$L:=\Phi(-q^9)^2 + 2 q\Phi(-q^9) X(-q^3) +
 4 q^2 X(-q^3)^2 $$ and $$ M:=P(q^3)^2 - q\Psi(q^9) P(q^3) + q^2\Psi(q^9)^2.$$
 Then all terms having the exponents of the form of $3n+2 (n\geq 0)$ in powers of $q$ in $L^4M^4$ are
 $A+B+C+D+E$, Where
\begin{eqnarray*}
 A&=& 40q^2P(q^3)^8 X(-q^3)^2\Phi(-q^9)^6-32q^2P(q^3)^7 X(-q^3)\Phi(-q^9)^7\Psi(q^9)\cr
 &+&10q^2P(q^3)^6\Phi(-q^9)^8\Psi(q^9)^2,\cr
 B&=& 512q^5P(q^3)^8 X(-q^3)^5\Phi(-q^9)^3 -1216q^5 P(q^3)^7 X(-q^3)^4\Phi(-q^9)^4\Psi(q^9)\cr
 &+&1280q^5P(q^3)^6 X(-q^3)^3\Phi(-q^9)^5\Psi(q^9)^2\cr
 &-&640q^5 P(q^3)^5 X(-q^3)^2\Phi(-q^9)^6\Psi(q^9)^3\cr
 &+& 152q^5P(q^3)^4 X(-q^3)^3\Phi(-q^9)^7\Psi(q^9)^4-16q^5P(q^3)^3 \Phi(-q^9)^8\Psi(q^9)^5,\cr
C&=&256q^8P(q^3)^8X(-q^3)^8-2048q^8P(q^3)^7 X(-q^3)^7 \Phi(-q^9)\Psi(q^9) \cr
 &+& 6400q^8P(q^3)^6 X(-q^3)^6 \Phi(-q^9)^2\Psi(q^9)^2\cr
 &-&8192q^8P(q^3)^5 X(-q^3)^5 \Phi(-q^9)^3\Psi(q^9)^3\cr
&+& 5776q^8P(q^3)^4 X(-q^3)^4 \Phi(-q^9)^4\Psi(q^9)^4\cr
&-&2048q^8P(q^3)^3 X(-q^3)^3 \Phi(-q^9)^5\Psi(q^9)^5\cr
 &+& 400q^8P(q^3)^2 X(-q^3)^2 \Phi(-q^9)^6\Psi(q^9)^6-32q^8P(q^3) X(-q^3) \Phi(-q^9)^7\Psi(q^9)^7\cr
 &+&q^8 \Phi(-q^9)^8\Psi(q^9)^8,\cr
 D&=&-4096q^{11}P(q^3)^5 X(-q^3)^8 \Psi(q^9)^3+9728q^{11}P(q^3)^4 X(-q^3)^7 \Phi(-q^9)\Psi(q^9)^4 \cr
 &-&10240q^{11}P(q^3)^3 X(-q^3)^6 \Phi(-q^9)^2\Psi(q^9)^5\cr
 &+&5120q^{11}P(q^3)^2 X(-q^3)^5 \Phi(-q^9)^3\Psi(q^9)^6\cr
 & -&1216q^{11}P(q^3) X(-q^3)^4 \Phi(-q^9)^4\Psi(q^9)^7+128q^{11} X(-q^3)^3 \Phi(-q^9)^5\Psi(q^9)^8
  \end{eqnarray*} and
 \begin{eqnarray*}
 E&=& 2560q^{14}P(q^3)^2 X(-q^3)^8 \Psi(q^9)^6-2048q^{14}P(q^3) X(-q^3)^7 \Phi(-q^9)\Psi(q^9)^7\cr
 &+& 640q^{14} X(-q^3)^6 \Phi(-q^9)^2\Psi(q^9)^8.
\end{eqnarray*}
\end{lemma}
\begin{proof}
We directly expand the expression of $L^4M^4$ and then extract the terms having exponents $3n+2$ in powers of $q$, then
we can obtain the results above.
\end{proof}
\begin{lemma}\label{4.2}
\begin{eqnarray}
A&=& 2\cdot 3^2 q^2\frac{(q^6;q^6)_\infty^6(q^9;q^9)_\infty^{26}}{(q^3;q^3)_\infty^6(q^{18};q^{18})_\infty^{10}},\cr
B&=& 8\cdot 3^2 q^5 \frac{(q^6;q^6)_\infty^3(q^9;q^9)_\infty^{17}}{(q^3;q^3)_\infty^3(q^{18};q^{18})_\infty},\cr
C&=& 19 \cdot 3^3 q^8 (q^9;q^9)_\infty^8(q^{18};q^{18})_\infty^{8},\cr
D&=& -64 \cdot 3^2 q^{11} \frac{(q^3;q^3)_\infty^3(q^{18};q^{18})_\infty^{17}}{(q^6;q^6)_\infty^3(q^{9};q^{9})_\infty},\cr
E&=& 128\cdot 3^2 q^{14} \frac{(q^3;q^3)_\infty^6(q^{18};q^{18})_\infty^{26}}{(q^6;q^6)_\infty^6(q^{9};q^{9})_\infty^{10}}.
\end{eqnarray}
\end{lemma}
\begin{proof}
By using the definitions of $P(q)$, $X(-q)$, the Lemma \ref{lemma2.2}, Lemma \ref{lemma2.4} and the formulas $$(q^6;q^6)_\infty=(q^6;q^{18})_\infty(q^{12};q^{18})_\infty(q^{18};q^{18})_\infty,\quad (q^9;q^9)_\infty=(q^{9};q^{18})_\infty(q^{18};q^{18}),$$ we have
\begin{eqnarray*}
A&=& 40q^2\frac{(q^6;q^{18})_\infty^{8}(q^{12};q^{18})_\infty^8(q^9;q^9)_\infty^{16}(q^3;q^3)_\infty^2(q^{18};q^{18})_\infty^4(q^9;q^9)_\infty^{12}}
{(q^3;q^3)_\infty^8(q^6;q^6)_\infty^2(q^9;q^9)_\infty^2(q^{18};q^{18})_\infty^6}\cr
&-& 32q^2 \frac{(q^6;q^{18})_\infty^{7}(q^{12};q^{18})_\infty^7(q^9;q^9)_\infty^{14}(q^3;q^3)_\infty
(q^{18};q^{18})_\infty^2(q^9;q^9)_\infty^{14}(q^{18};q^{18})_\infty}
{(q^3;q^3)_\infty^7(q^6;q^6)_\infty(q^9;q^9)_\infty(q^{18};q^{18})_\infty^7(q^9;q^{18})_{\infty}}\cr
&+&10q^2\frac{(q^6;q^{18})_\infty^{6}(q^{12};q^{18})_\infty^6(q^9;q^9)_\infty^{12}(q^9;q^9)_\infty^{16}(q^{18};q^{18})_\infty^2}
{(q^3;q^3)_\infty^6(q^{18};q^{18})_\infty^8(q^9;q^{18})_\infty^2}\cr
&=& 18 q^2 \frac{(q^6;q^6)_\infty^6(q^9;q^9)_\infty^{26}}{(q^3;q^3)_\infty^6(q^{18};q^{18})_\infty^{10}}.
\end{eqnarray*}
\begin{eqnarray*}
B&=& 512q^5 \frac{(q^6;q^{18})_\infty^{8}(q^{12};q^{18})_\infty^8(q^9;q^9)_\infty^{16}(q^3;q^3)_\infty^5(q^{18};q^{18})_\infty^{10}(q^9;q^9)_\infty^6}
{(q^3;q^3)_\infty^8(q^6;q^6)_\infty^5(q^9;q^9)_\infty^5(q^{18};q^{18})_\infty^3}\cr
&-& 1216q^5\frac{(q^6;q^{18})_\infty^{7}(q^{12};q^{18})_\infty^7(q^9;q^9)_\infty^{14}(q^3;q^3)_\infty^4(q^{18};q^{18})_\infty^{8}(q^9;q^9)_\infty^8
(q^{18};q^{18})_\infty}
{(q^3;q^3)_\infty^7(q^6;q^6)_\infty^4(q^9;q^9)_\infty^4(q^{18};q^{18})_\infty^4(q^9;q^{18})_\infty}\cr
&+& 1280q^5 \frac{(q^6;q^{18})_\infty^{6}(q^{12};q^{18})_\infty^6(q^9;q^9)_\infty^{12}(q^3;q^3)_\infty^3(q^{18};q^{18})_\infty^{6}
(q^9;q^9)_\infty^{10}(q^{18};q^{18})_\infty^2}
{(q^3;q^3)_\infty^6(q^6;q^6)_\infty^3(q^9;q^9)_\infty^3(q^{18};q^{18})_\infty^5(q^9;q^{18})_\infty^2}\cr
&-& 640q^5 \frac{(q^6;q^{18})_\infty^{5}(q^{12};q^{18})_\infty^5(q^9;q^9)_\infty^{10}(q^3;q^3)_\infty^2(q^{18};q^{18})_\infty^{4}
(q^9;q^9)_\infty^{12}(q^{18};q^{18})_\infty^3}
{(q^3;q^3)_\infty^5(q^6;q^6)_\infty^2(q^9;q^9)_\infty^2(q^{18};q^{18})_\infty^6(q^9;q^{18})_\infty^3}\cr
&+& 152 q^5 \frac{(q^6;q^{18})_\infty^{4}(q^{12};q^{18})_\infty^4(q^9;q^9)_\infty^{8}(q^3;q^3)_\infty(q^{18};q^{18})_\infty^{2}
(q^9;q^9)_\infty^{14}(q^{18};q^{18})_\infty^4}
{(q^3;q^3)_\infty^4(q^6;q^6)_\infty(q^9;q^9)_\infty(q^{18};q^{18})_\infty^7(q^{9};q^{18})_\infty^4}\cr
&-& 16q^5\frac{(q^6;q^{18})_\infty^{3}(q^{12};q^{18})_\infty^3(q^9;q^9)_\infty^{6}(q^9;q^9)_\infty^{16}(q^{18};q^{18})_\infty^5}
{(q^3;q^3)_\infty^3(q^{18};q^{18})_\infty^6(q^{9};q^{18})_\infty^5}\cr
&=& 8\cdot 3^2  q^5 \frac{(q^6;q^6)_\infty^3(q^9;q^9)_\infty^{17}}{(q^3;q^3)_\infty^3(q^{18};q^{18})_\infty}.
\end{eqnarray*}
\begin{eqnarray*}
C&=& 256 q^8 \frac{(q^6;q^{18})_\infty^{8}(q^{12};q^{18})_\infty^8(q^9;q^9)_\infty^{16}(q^3;q^3)_\infty^8(q^{18};q^{18})_\infty^{16}}
{(q^3;q^3)_\infty^8(q^6;q^6)_\infty^8(q^9;q^9)_\infty^8} \cr
    &-&2048 q^8 \frac{(q^6;q^{18})_\infty^{7}(q^{12};q^{18})_\infty^7(q^9;q^9)_\infty^{14}(q^3;q^3)_\infty^7(q^{18};q^{18})_\infty^{14}
(q^9;q^9)_\infty^{2}(q^{18};q^{18})_\infty}
{(q^3;q^3)_\infty^7(q^6;q^6)_\infty^7(q^9;q^9)_\infty^7(q^{18};q^{18})_\infty(q^{9};q^{18})_\infty} \cr
    &+&6400 q^8 \frac{(q^6;q^{18})_\infty^{6}(q^{12};q^{18})_\infty^6(q^9;q^9)_\infty^{12}(q^3;q^3)_\infty^6(q^{18};q^{18})_\infty^{12}
(q^9;q^9)_\infty^{4}(q^{18};q^{18})_\infty^2}
{(q^3;q^3)_\infty^6(q^6;q^6)_\infty^6(q^9;q^9)_\infty^6(q^{18};q^{18})_\infty^2(q^{9};q^{18})_\infty^2}  \cr
   &-& 8192 q^8 \frac{(q^6;q^{18})_\infty^{5}(q^{12};q^{18})_\infty^5(q^9;q^9)_\infty^{10}(q^3;q^3)_\infty^5(q^{18};q^{18})_\infty^{10}
 (q^9;q^9)_\infty^{6}(q^{18};q^{18})_\infty^3}
 {(q^3;q^3)_\infty^5(q^6;q^6)_\infty^5(q^9;q^9)_\infty^5(q^{18};q^{18})_\infty^3(q^{9};q^{18})_\infty^3}\cr
   &+&5776 q^8\frac{(q^6;q^{18})_\infty^{4}(q^{12};q^{18})_\infty^4(q^9;q^9)_\infty^{8}(q^3;q^3)_\infty^4
(q^{18};q^{18})_\infty^{8}(q^9;q^9)_\infty^{8}(q^{18};q^{18})_\infty^4}
{(q^3;q^3)_\infty^4(q^6;q^6)_\infty^4(q^9;q^9)_\infty^4(q^{18};q^{18})_\infty^4(q^{9};q^{18})_\infty^4}\cr
   &-&2048 q^8\frac{(q^6;q^{18})_\infty^{3}(q^{12};q^{18})_\infty^3(q^9;q^9)_\infty^{6}(q^3;q^3)_\infty^3(q^{18};q^{18})_\infty^{6}
(q^9;q^9)_\infty^{10}(q^{18};q^{18})_\infty^5}
{(q^3;q^3)_\infty^3(q^6;q^6)_\infty^3(q^9;q^9)_\infty^3(q^{18};q^{18})_\infty^5(q^{9};q^{18})_\infty^5}\cr
   &+&400q^8 \frac{(q^6;q^{18})_\infty^{2}(q^{12};q^{18})_\infty^2(q^9;q^9)_\infty^{4}(q^3;q^3)_\infty^2(q^{18};q^{18})_\infty^{4}
(q^9;q^9)_\infty^{12}(q^{18};q^{18})_\infty^6}
{(q^3;q^3)_\infty^2(q^6;q^6)_\infty^2(q^9;q^9)_\infty^2(q^{18};q^{18})_\infty^6(q^{9};q^{18})_\infty^6}\cr
  &-&32q^8  \frac{(q^6;q^{18})_\infty^{}(q^{12};q^{18})_\infty(q^9;q^9)_\infty^{2}(q^3;q^3)_\infty(q^{18};q^{18})_\infty^{2}
 (q^9;q^9)_\infty^{14}(q^{18};q^{18})_\infty^7}
 {(q^3;q^3)_\infty(q^6;q^6)_\infty(q^9;q^9)_\infty(q^{18};q^{18})_\infty^7(q^{9};q^{18})_\infty^7}\cr
  &+& q^8\frac{(q^9;q^9)_\infty^{16}(q^{18};q^{18})_\infty^8}{(q^{18};q^{18})_\infty^8(q^{9};q^{18})_\infty^8}\cr
  &=& 19\cdot 3^3 q^8(q^9;q^9)_\infty^{8}(q^{18};q^{18})_\infty^8.\cr\cr
D&=& -4096 q^{11}\frac{(q^6;q^{18})_\infty^{5}(q^{12};q^{18})_\infty^5(q^9;q^9)_\infty^{10}(q^3;q^3)_\infty^8(q^{18};q^{18})_\infty^{16}
(q^{18};q^{18})_\infty^3}{(q^3;q^3)_\infty^5(q^6;q^6)_\infty^8(q^9;q^9)_\infty^8(q^{9};q^{18})_\infty^3}\cr
 &+&9728q^{11} \frac{(q^6;q^{18})_\infty^{4}(q^{12};q^{18})_\infty^4(q^9;q^9)_\infty^{8}(q^3;q^3)_\infty^7(q^{18};q^{18})_\infty^{14}
(q^9;q^9)_\infty^{2}(q^{18};q^{18})_\infty^4}{(q^3;q^3)_\infty^4(q^6;q^6)_\infty^7(q^9;q^9)_\infty^7(q^{18};q^{18})_\infty(q^{9};q^{18})_\infty^4}\cr
 &-&10240q^{11}\frac{(q^6;q^{18})_\infty^{3}(q^{12};q^{18})_\infty^3(q^9;q^9)_\infty^{6}(q^3;q^3)_\infty^6(q^{18};q^{18})_\infty^{12}
(q^9;q^9)_\infty^{4}(q^{18};q^{18})_\infty^5}{(q^3;q^3)_\infty^3(q^6;q^6)_\infty^6(q^9;q^9)_\infty^6(q^{18};q^{18})_\infty^2(q^{9};q^{18})_\infty^5}
               \\
 &+&5120q^{11}\frac{(q^6;q^{18})_\infty^{2}(q^{12};q^{18})_\infty^2(q^9;q^9)_\infty^{4}(q^3;q^3)_\infty^5(q^{18};q^{18})_\infty^{10}
(q^9;q^9)_\infty^{6}(q^{18};q^{18})_\infty^6}
{(q^3;q^3)_\infty^2(q^6;q^6)_\infty^5(q^9;q^9)_\infty^5(q^{18};q^{18})_\infty^3(q^{9};q^{18})_\infty^6}
 \end{eqnarray*}
 \begin{eqnarray*}
 &-&1216q^{11}\frac{(q^6;q^{18})_\infty^{}(q^{12};q^{18})_\infty(q^9;q^9)_\infty^{2}(q^3;q^3)_\infty^4(q^{18};q^{18})_\infty^{8}
(q^9;q^9)_\infty^{8}(q^{18};q^{18})_\infty^7}
{(q^3;q^3)_\infty(q^6;q^6)_\infty^4(q^9;q^9)_\infty^4(q^{18};q^{18})_\infty^4(q^{9};q^{18})_\infty^7}\cr
 &+&128q^{11}\frac{(q^3;q^3)_\infty^3(q^{18};q^{18})_\infty^{6}(q^9;q^9)_\infty^{10}(q^{18};q^{18})_\infty^8}
{(q^6;q^6)_\infty^3(q^9;q^9)_\infty^3(q^{18};q^{18})_\infty^5(q^{9};q^{18})_\infty^8}\cr
 &=& -64\cdot 3^2 q^{11} \frac{(q^3;q^3)_\infty^{3}(q^{18};q^{18})_\infty^{17}}{(q^{6};q^{6})_\infty^3(q^{9};q^{9})_\infty}
\end{eqnarray*}
\begin{eqnarray*}
E&=& 2560q^{14}\frac{(q^6;q^{18})_\infty^{2}(q^{12};q^{18})_\infty^2(q^9;q^9)_\infty^{4}(q^3;q^3)_\infty^8(q^{18};q^{18})_\infty^{16}
(q^{18};q^{18})_\infty^6}{(q^3;q^3)_\infty^2(q^6;q^6)_\infty^8(q^9;q^9)_\infty^8(q^{9};q^{18})_\infty^6}\cr
&-& 2048 q^{14}\frac{(q^6;q^{18})_\infty^{}(q^{12};q^{18})_\infty(q^9;q^9)_\infty^{2}(q^3;q^3)_\infty^7
(q^{18};q^{18})_\infty^{14}(q^9;q^9)_\infty^{2}(q^{18};q^{18})_\infty^7}
{(q^3;q^3)_\infty(q^6;q^6)_\infty^7(q^9;q^9)_\infty^7(q^{18};q^{18})_\infty(q^{9};q^{18})_\infty^7}\cr
&+& 640 q^{14}\frac{(q^3;q^3)_\infty^6(q^{18};q^{18})_\infty^{12}(q^9;q^9)_\infty^{4}(q^{18};q^{18})_\infty^8}
{(q^6;q^6)_\infty^6(q^9;q^9)_\infty^6(q^{18};q^{18})_\infty^2(q^{9};q^{18})_\infty^8}\cr
&=& 128\cdot 3^2 q^{14} \frac{(q^3;q^3)_\infty^6(q^{18};q^{18})_\infty^{26}}{(q^6;q^6)_\infty^6(q^9;q^9)_\infty^{10}}.
\end{eqnarray*}
\end{proof}

Now we prove the identity \ref{xiong}.
\begin{proof}By Theorem \ref{chanxizhi}, Lemma \ref{4.1} and Lemma \ref{4.2}, we have
\begin{eqnarray*}
\sum_{n=0}^\infty a(9n+8)q^{3n+2}&=&\text{ all terms having the exponents of the form of}\,\, 3n+2\\
 &&(n\geq 0)\,\, \text{in powers of}\,\, q\,\, \text{in}\,\,3\frac{(q^3;q^3)_\infty^3(q^6;q^6)_\infty^3}{(q;q)_{\infty}^4(q^2;q^2)_{\infty}^4}\cr\cr
&=& \text{ all terms having the exponents of the form of}\,\, 3n+2\\&& (n\geq 0)\,\, \text{in powers of}\,\, q\,\, \text{in}\,\,3\frac{(q^3;q^3)_\infty^3(q^6;q^6)_\infty^3}{\Phi(-q)^4\Psi(q)^4}\cr\cr
&=& \text{ all terms having the exponents of the form of}\,\, 3n+2\\&& (n\geq 0)\,\, \text{in powers of}\,\, q\,\, \text{in}\,\,3(q^3;q^3)_\infty^3(q^6;q^6)_\infty^3\\&&\cdot\frac{\Phi(-q^9)^4\Psi(q^9)^4}{\Phi(-q^3)^{16}\Psi(q^3)^{16}}\cdot L^4M^4\cr\cr
&=& 3(q^3;q^3)_\infty^3(q^6;q^6)_\infty^3\cdot\frac{\Phi(-q^9)^4\Psi(q^9)^4}{\Phi(-q^3)^{16}\Psi(q^3)^{16}}\\&&\cdot
(A+B+C+D+E)\cr\cr
&=& 3\frac{(q^9;q^9)_\infty^{4}(q^{18};q^{18})_\infty^4}{(q^3;q^3)_\infty^{16}(q^6;q^6)_\infty^{16}}\cdot (A+B+C+D+E)\cr\cr
&=& 3\frac{(q^9;q^9)_\infty^{4}(q^{18};q^{18})_\infty^4}{(q^3;q^3)_\infty^{16}(q^6;q^6)_\infty^{16}}
\cdot (18 q^2 \frac{(q^6;q^6)^6(q^9;q^9)^{26}}{(q^3;q^3)^6(q^{18};q^{18})^{10}}\cr\cr&+&
8\cdot 3^2  q^5 \frac{(q^6;q^6)_\infty^3(q^9;q^9)_\infty^{17}}{(q^3;q^3)_\infty^3(q^{18};q^{18})_\infty}+
19\cdot 3^3 q^8(q^9;q^9)_\infty^{8}(q^{18};q^{18})_\infty^8\cr\cr
&-&64\cdot 3^2 q^{11} \frac{(q^3;q^3)_\infty^{3}(q^{18};q^{18})_\infty^{17}}{(q^{6};q^{6})_\infty^3(q^{9};q^{9})_\infty}\\
&+&128\cdot 3^2 q^{14} \frac{(q^3;q^3)_\infty^6(q^{18};q^{18})_\infty^{26}}{(q^6;q^6)_\infty^6(q^9;q^9)_\infty^{10}}
)
\end{eqnarray*}
If we divide $q^2$ on both sides and replace $q^3$ by $q$, we obtain
\begin{eqnarray*}
\sum_{n=0}^{\infty}a(9n+8)q^n &=& 2\cdot 3^3 \frac{(q^3;q^3)_{\infty}^{30}}{(q;q)_{\infty}^{19}(q^2;q^2)_{\infty}^{7}(q^6;q^6)_{\infty}^{6}}+8\cdot 3^3 q \frac{(q^3;q^3)_{\infty}^{21}(q^6;q^6)_{\infty}^{3}}{(q;q)_{\infty}^{16}(q^2;q^2)_{\infty}^{10}}\cr\cr&+&19\cdot 3^4 q^2 \frac{(q^3;q^3)_{\infty}^{12}(q^6;q^6)_{\infty}^{12}}{(q;q)_{\infty}^{13}(q^2;q^2)_{\infty}^{13}}
-64\cdot 3^{3} q^3 \frac{(q^3;q^3)_{\infty}^{3}(q^6;q^6)_{\infty}^{21}}{(q;q)_{\infty}^{10}(q^2;q^2)_{\infty}^{16}}\cr\cr&+&128\cdot 3^{3}q^4 \frac{(q^6;q^6)_{\infty}^{30}}{(q;q)_{\infty}^{7}(q^2;q^2)_{\infty}^{19}(q^3;q^3)_{\infty}^{6}},
\end{eqnarray*}
which is the identity \ref{xiong}.
\end{proof}

\section{ACKNOWLEDGEMENT}
The author thanks Professor Liu Zhi-Guo's helpful comments.

\bibliographystyle{amsplain}

\end{document}